\documentclass{amsart}
\usepackage{amssymb}
\usepackage{amsmath}
\usepackage{amsfonts}
\usepackage{graphicx}
\usepackage{hyperref}
\usepackage{fancyhdr}
\makeatletter
\providecommand{\U}[1]{\protect\rule{.1in}{.1in}}
\newtheorem{theorem}{Theorem}

\newtheorem{remark}[theorem]{Remark}

\begin{document}
\date{\today}
\title{Huygens synchronization of three clocks equidistant from each other}
\author{Emma D'Aniello$^1$ and Henrique M. Oliveira$^{2*}$}
\thanks{$^*$ Corresponding author}
\thanks{$^1$ORCID: 0000-0001-5872-0869 Dipartimento di Matematica e Fisica, Universit\`{a} degli Studi della Campania
\textquotedblleft Luigi Vanvitelli\textquotedblright, Viale Lincoln n. 5 -
81100 Caserta, Italia}

\thanks{$^2$ORCID: 0000-0002-3346-4915
Department of Mathematics \and Center for Mathematical Analysis, Geometry and
Dynamical Systems, Instituto Superior T\'ecnico, University of Lisbon, Av.
Rovisco Pais, 1049-001, Lisboa, Portugal }
\thanks{$^1$emma.daniello@unicampania.it; $^2$holiv@math.tecnico.ulisboa.pt}
\thanks{The author ED was partially supported by the program Erasmus+. The author HMO was partially supported by  Fundação para a Ciência e
Tecnologia, UIDB/04459/2020 and UIDP/04459/2020.
}

\begin{abstract}
This paper investigates the synchronization of three identical oscillators, or clocks, suspended from a common rigid support. We consider scenarios where each clock interacts with the other two, achieving synchronization through small impacts exchanged between oscillator pairs. The fundamental outcome of our study reveals that the ultimate synchronized state maintains a phase difference of $\frac{2\pi}{3}$ between successive clocks, either clockwise or counter-clockwise. Furthermore, these locked states exhibit an attracting set, which closure encompasses the entire initial conditions space. Our analytical approach involves constructing a nonlinear discrete dynamical system in dimension two.

These findings hold significance for sets of three weakly coupled periodic oscillators engaged in mutual symmetric impact periodic interaction, irrespective of the specific oscillator models employed. Lastly, we explore the amplitude of oscillations at the final locked state in the context of two and three interacting Andronov pendulum clocks. Our analysis reveals a precise small increase in the amplitude of the locked-state oscillations, as quantified in this paper.
\end{abstract}
\maketitle

\section{Introduction}\label{intro}

Synchronization among oscillators with some form of coupling is commonly referred to as universal \cite{strogatz2004} and plays a significant role in natural phenomena \cite{carvalho2020,Luo2009,Luo2013,Pit}.

In 1665, Christiaan Huygens, the inventor of the pendulum clock, made a noteworthy observation regarding synchronization between two pendulum clocks suspended from the same support \cite{Huy}. While confined to bed due to illness, he witnessed both in-phase and anti-phase synchronization as the final states of the coupled system. Huygens later replicated this experiment by suspending the two clocks on a board supported by chairs.

These distinct observations by Huygens paved the way for two separate lines of analysis as is acknowledged in \cite{goldsztein2021}. The subsequent exploration involved clocks attached to a wall with momentum conservation in the clocks-beam system, allowing for the movement of the supporting plank. This line of study culminated in numerous investigations \cite{Benn,Col,Frad, Jova,Col2,Martens,Oud,Ram2,Sen}. We term this model the "classical model," deeply rooted in classical mechanics and accounting for viscous friction.

In this paper, we depart from the classical model inspired by the works of Vassalo-Pereira \cite{Vass} and Ralf Abraham and co-authors \cite{Abr2,Abr}, utilizing Andronov's established model for the pendulum clock with dry friction.

Our focus is the synchronization of oscillators in a plane with an asymptotically stable limit cycle. We specifically consider isochronous clocks – those with a frequency independent of amplitude. We investigate the scenario where pendulums are suspended from a rigid house beam, rendering them immobile. To approach this, we employ a perturbative model. This model accounts for the interaction between the pendulums arising from their internal impacts. The ensuing perturbation generates traveling waves that transfer kinetic energy between the oscillators.

Our study takes inspiration from a model proposed in \cite{OlMe}. This model introduced a theoretical framework for such interactions, accompanied by simulations and experiments involving an immobile support wall. Unlike the classical model, our approach involves discrete impacts exchanged between three identical oscillators. These impacts occur once per cycle, in the form of short burst travelling waves. We assume instantaneous coupling due to the rapid speed of mechanical waves in the medium to which the clocks are attached \cite{OlMe}.

Our investigation reveals a symmetric asymptotic state where all clocks maintain a phase difference of $\frac{2\pi}{3}$ between each other. We term this phenomenon "Huygens synchronization," extending the concept from \cite{Abr2,Abr,Vass,OlMe}. 

The structure of this paper comprises five sections. Section \ref{sec2} delves into the original pendulum clock model, briefly recapping the model for two identical clocks. In Section \ref{sec3}, we derive the model for three identical clocks with mutual interactions. Detailed construction information is provided in the appendix. Section \ref{sec4} encompasses an analysis of the model's symmetries and stabilities. Finally, Section \ref{sec5} offers conclusions and outlines potential directions for future research.

\section{Model for the synchronization of two oscillators}\label{sec2}

\subsection{Some background}

For the sake of completeness, we present a concise theory of synchronization for two oscillators exchanging small perturbations at each cycle. We focus on identical oscillators, and this theory's applicability extends to networks of identical oscillators, electronic oscillators, and various other real-world systems. In future investigations, we aim to explore cases involving slightly different oscillators, which lead to regions of stability versus instability in the parameter space, known as Arnold Tongues \cite{boyland1986,gilmore2011,OlPe2022}.

For fundamental and classical definitions and concepts related to synchronization, such as phase and frequency, we follow and refer to \cite{Pit}. For broader concepts concerning the general theory of dynamical systems, such as limit cycles, we refer to \cite{arrowsmith1990}. Throughout this paper, we consider oscillators as dynamical systems exhibiting limit cycles. We use the term "clock" to refer to a specific type of oscillator as described by the Andronov model \cite{OlMe}.

Given a point $p_{0}$ on the limit cycle $\gamma$, the time required for a return to $p_{0}$ after completing one cycle on the limit cycle is denoted by the period $T_{0}$. A phase $\varphi$ serves as a real coordinate describing the representative point's position on the limit cycle \cite{Nakao,Pit}.

Let $B_{\gamma}$ represent the basin of attraction of the limit cycle. For points outside the limit cycle $\gamma$ but within $B_{\gamma}$, we extend the phase definition as follows: all points $p$ in $B_{\gamma}$ that converge to the same $p_{0}$ on the limit cycle $\gamma$ as $t\rightarrow\infty$ are assigned the phase $\varphi$ of $p_{0}$ \cite{gu1975}. The set of points sharing the same phase forms an isochron curve. When oscillator states lie on the same isochron at a given time, they remain on the same isochron over time \cite{gu1975,Nakao}. In the presence of perturbations, each clock's state can deviate slightly from the limit cycle and generally jump to another isochron. We also assume that the limit cycles remain structurally stable under minor perturbations.

Considering two oscillators, labeled $1$ and $2$, with orbits on or near the limit cycle, each has a distinct phase, denoted $\varphi$ and $\psi$, respectively.

Studying the synchronization of these oscillators involves establishing a dynamical system for the phase difference between them.

There are two potential research directions \cite{Pit}. The first examines phase differences over continuous time, expressed as the function $\phi\left( t\right) =\psi\left( t\right) -\varphi\left( t\right) $\ for $t\in\left[ 0,+\infty\right[ $. The second, which we adopt in this paper, analyzes the discrete phase difference $\phi_{n}=\psi_{n} -\varphi_{n}$ at specific instances $n=0,1,2,\ldots$. This study exclusively employs the latter approach.

Phase synchronization occurs when the phase differences between oscillators converge toward a distinct attractor. When this attractor is a solitary point, phase locking is established. Naturally, more intricate coupled states can arise \cite{Martens}. The primary objective of any synchronization theory is to derive the dynamics of this phase difference and establish the nature and existence of the attractor. In the context of Huygens' observations, the attractor was either the point $0$ or the point $\pi$, and the phase dynamics remained unidimensional.

\subsection{The Andronov model for an isolated clock}

We revisit the model that assumes the prevalence of dry friction within the internal metal components of the clock, with viscous damping playing a secondary role. Using the angular coordinate $q$, the differential equation governing the isolated pendulum clock is given by
\begin{equation}
\ddot{q}+\mu\text{ }\operatorname*{sign}\dot{q}+q=0,
\end{equation}
where $\mu>0$ represents the dry friction coefficient, and $\operatorname*{sign}\left(
x\right)  $ is the classical function that takes the value $-1$ for $x<0$ and $1$
for $x>0$. In \cite{And}, it was considered that during each cycle, the escape
mechanism imparts a fixed amount of normalized kinetic energy $\frac{h^{2}}{2}$ to the pendulum, thus compensating for the kinetic energy loss caused by dry friction in each complete cycle. This transfer of kinetic energy is termed a "kick." The origin is positioned so that the kick is given precisely when $q=-\mu$. The phase portrait is depicted in Fig. $1$.

\begin{figure}[ht]
\begin{center}
\includegraphics[
height=2.1075in,
width=2.2753in
]%
{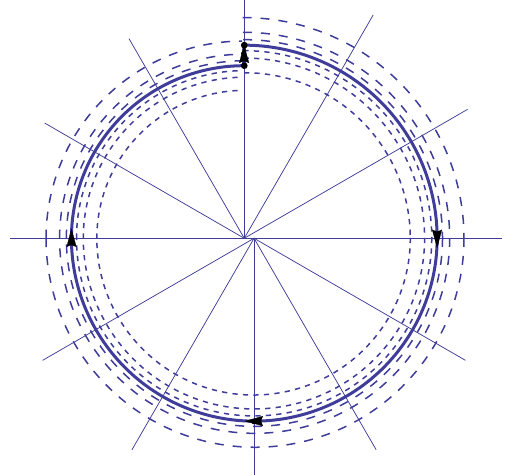}%
\caption{Limit cycle of an isolated clock represented as a solid curve in the
phase space. The horizontal axis represents the angular position, and the
vertical axis represents velocity.}%
\label{Fig1}%
\end{center}
\end{figure}

Similar to \cite{OlMe}, with initial conditions $q\left(  t=0\right)  =-\mu$ and
$\dot{q}\left(  t=0\right)  =v_{0}$, a Poincaré section (see vol. II, page
268, of \cite{Bir}) is defined as the half line $q=-\mu^{+}$ and $\dot
{q}>0$ \cite{And}. The symbol $+$ indicates that we are considering the section immediately after the kick. Due to friction over a complete cycle, a velocity reduction of $-4\mu$ takes place. By evaluating the velocity,
$v_{n}=\dot{q}\left(  2n\pi^{+}\right)  $, at the Poincaré section in each
cycle, the non-linear discrete dynamical system \cite{And} is derived as
\begin{equation}
v_{n+1}=\sqrt{\left(  v_{n}-4\mu\right)  ^{2}+h^{2}}\text{.}\label{recu}%
\end{equation}
This equation possesses the asymptotically stable fixed point
\begin{equation}
v_{f}=\frac{h^{2}}{8\mu}+2\mu\text{.}%
\end{equation}
Any initial condition $v_{0}\in\left(  4\mu,+\infty\right)  $ converges to
$v_{f}$. Each cycle corresponds to a phase increment of $2\pi$, and the phase
$\varphi$ linearly varies with respect to $t$, specifically
\[
\varphi=2\pi t.
\]
As already mentioned, the nature of limit cycle is not of fundamental
importance when we consider the interaction of three identical clocks, as we
shall see in the sequel. We have presented here the basis of our reasoning in
the non-usual case when the computations of the limit cycle are explicit and
the usual angular phase is a linear function of $t$. 

We define the amplitude of the movement of the pendulum for the Andronov model as exactly
this value $v_{f}$, which is the maximum of the angular velocity and
the maximum of the angle coordinate $q$.

\subsection{Two interacting oscillators}

We briefly present the model constructed in \cite{OlMe} for two pendulums with
the same natural frequency. When one clock receives the kick, the impact
propagates in the wall slightly perturbing the second clock. The perturbation
is assumed to be instantaneous since the time of travel of sound, 
i.e., the mechanical waves as explained in the introduction, in the wall between the clocks is 
assumed very small compared to
the period. This reasoning allows a treatment of the three close
clocks as if they were equidistant even if it is not exactly the case. 

To describe and investigate the effect of the kicks, we construct a discrete
dynamical system for the phase difference between the two clocks. We compute
each cycle using as reference one of the clocks (the choice is irrelevant,
since the model is symmetric). We choose, to fix ideas, clock 1 as the
reference: whenever its phase reaches $0$ $\left(  \operatorname{mod}%
2\pi\right)  $, the number of cycles increases one unit from $n$ to $n+1$.

If there exists an attracting fixed point for that dynamical system, the phase
locking occurs.

The secular repetition of perturbations leads the system with the two clocks
in phase opposition as Huygens observed in 1665 \cite{Huy}.

In the case of frequencies equal to one and small friction coefficient, the
discrete dynamical model obtained in \cite{OlMe} for the phase difference
between two clocks, $\phi_{n}=\psi_{n}-\varphi_{n}$, gives the Adler equation
\cite{adler1946study,Pit}%
\begin{equation}
\phi_{n+1}=\phi_{n}+\varepsilon\sin\phi_{n}, \label{Perturb}%
\end{equation}
with a very small constant $\varepsilon=\frac{16\mu\alpha}{h^{2}}$, where
$\alpha$ is a small coupling constant measuring the mutual perturbations. In
the interval $\left[  0,2\pi\right[  $, there are two fixed points which are
$\pi$ and $0$, attracting and repelling, respectively.

Equation (\ref{Perturb}) is the starting point from where we begin, in the
present paper, the study the three symmetric clocks in mutual interaction.

\begin{remark}
In any model with a perturbation of phase given by equation (\ref{Perturb})
per cycle, i.e., Adler's perturbation \cite{adler1946study,Pit}, despite being
a physical clock (with Andronov model or any different model) or other type of
oscillator, electric, quantic, electronic or biological, the theory presented
here for three oscillators interacting by small periodic impacts will be exactly the same, 
with the same conclusions.
\end{remark}

The amplitude of the oscillation in the case of two clocks increases slightly in the final state.
We focus our analysis on the velocity at $q=-\mu^{+}$ and $\dot
{q}=v_{f}$, i.e., at the Poincar\'{e} section. Since the mutual
interaction affects each clock when they are in phase opposition, the equation
for the Poincar\' e map for each generic oscillator is now

\begin{equation}
v_{n+1}=\sqrt{\left(  v_{n}-4\mu-\alpha\right)  ^{2}+h^{2}}\text{.}%
\end{equation}
This equation has the asymptotically stable\ fixed point%
\begin{equation}
v^{\ast}=\frac{\left(  4\mu+\alpha\right)  ^{2}+h^{2}}{2\left(  4\mu
+\alpha\right)  }\text{.}\label{perturb2}%
\end{equation}
The value of $v^{\ast}$, which is, in fact, the maximum of the
amplitude of the oscillation for the mutual interacting oscillators, is
slightly greater than $v_{f}$ of the isolated clock. More precisely, approximating 
in first order in $\alpha$, and keeping in mind that $h$ is
small, we get
\begin{equation}
v^{\ast}=v_{f}+\left(  \frac{1}{2}-\frac{h^{2}}{32\mu^{2}}\right)
\alpha.\label{perturb3}
\end{equation}

\section{Model for three pendulum clocks placed in the three vertices of an equilateral triangle}\label{sec3}

\subsection{Hypotheses}

We consider three pendulum clocks suspended at the same wall, placed in the
three vertices of an equilateral triangle, say the vertices are $A$, $B$, and
$C$ and $B$ are the extreme points of the basis of the triangles.

\begin{figure}
[h]
\begin{center}
\includegraphics[
height=2.7743in,
width=3.6443in
]%
{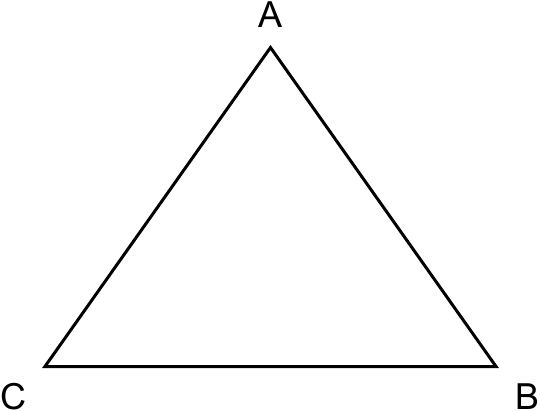}%
\caption{The three clocks hang at the three vertices of a triangle.}%
\label{Fig2}%
\end{center}
\end{figure}

This geometric setting is purely conceptual. Any set of three dynamical
systems receiving symmetric impacts from the other two will have the same type
of response of the clocks depicted in the three vertices of an equilateral triangle.

Call the clocks placed in the three vertices $A$, $B$ and $C$, respectively,
$O_{1}$, $O_{2}$ and $O_{3}$. When the clock $A$ receives the kick from the
escape mechanism, the impact propagates in the wall slightly perturbing the
other two clocks. As in \cite{OlMe}, the perturbation is assumed to be
instantaneous, since the time of travel of sound in the wall between the
clocks is assumed very small compared to the period. As for the two clocks
model discussed in \cite{OlMe}, we make the following assumptions, now
formulated for three clocks.

\begin{enumerate}
\item The system has dry friction \cite{And}.

\item \label{work3}The pendulums of clocks $O_{1},$ $O_{2}$ and $O_{3}$ have
respectively natural angular frequencies $\omega_{1}=\omega_{2}=\omega_{3}=1$.

\item The perturbation in the momentum is always in the same vertical
direction in the phase space \cite{Abr2, Abr}.

\item The friction coefficient is the same for all the three clocks, $\mu
_{1}=\mu_{2}=\mu_{3}=\mu$. The energy dissipated at each cycle of the three
clocks is the same, and the energy furnished by the escape mechanism to
compensate the loss of energy to friction in each cycle is $h_{1}=h_{2}%
=h_{3}=h$.

\item The perturbative interaction is instantaneous. This is a reasonable
assumption, since in general the perturbation propagation time between two
clocks is several orders of magnitude lower than the periods \cite{OlMe}.

\item \label{hyp6}The interaction is symmetric. The couplings have the same
constant $\alpha$ when one clock acts on another and conversely. In this model
$\alpha$ is assumed to be very small.

\item \label{hyp67}Each perturbation from clock $i$ to clock $j$ (where
$i,j\in\{1,2,3\}$ with $i\not =j$), when clock $i$ suffers its internal impact
of kinetic energy $h^{2}$, gives rise to a small perturbative change of phase
which is in first order a $2\pi$-periodic differentiable odd function $P$ of
the real variable $\phi$%
\begin{equation}
P\left(  \phi\right)  =\varepsilon\sin\phi\text{,} \label{Perturb1}%
\end{equation}
where $\phi=\phi_{ij}$ is the phase difference between clock $i$ and clock $j$.
\end{enumerate}

\begin{remark}
The value of the constant $\varepsilon$ in each interaction term is
$\varepsilon=\frac{8\mu\alpha}{h^{2}}$, i.e., half of the value obtained for
the two clocks \cite{OlMe}, where $\mu$ is the dry friction coefficient,
$\frac{h^{2}}{2}$ is the kinetic energy furnished by the internal escape
mechanism of each clock once per cycle and $\alpha$ is the interaction
coefficient between the clocks. The greater the $\alpha$ is, the greater the
mutual influence among the clocks. In this paper, we do not need to
particularize $\varepsilon$, since we are not interested in doing experimental
computations. In this paper, we are interested in the fundamental result of
symmetry between three oscillators subject to very weak mutual symmetric interaction.
\end{remark}

Most of the reasoning is independent on the form of the function $P\left(
\phi\right)  $, therefore we consider a general differentiable odd function of
the real variable $\phi$, $P(\phi)$, for the development of the model, and
consider it of the form (\ref{Perturb1}) when we analyse the model in section 4.

Observe that $|\sin(x+\varepsilon\sin y)-\sin x|<\varepsilon$ when
$\varepsilon$\ is assumed to be sufficiently small. Therefore, we restrict our
model to first order. We consider all the values of variables and constants in
IS units.

\subsection{Construction of the model}

We now construct a dynamical system using as reference the phase of the clock
in the vertex A (= clock $O_{1}$). This reference is arbitrary: any of the
clocks can be used as the reference clock with the same results at the end,
since the system is symmetric. We compute the effects of all phase differences
and perturbations when the clock at A makes a complete cycle returning to the
initial position. Without loss of generality, we consider the next working hypotheses.

\begin{enumerate}
\item \label{work1}The initial phase of clock at A at $t=0^{-}$ is zero, i.e.,
$\psi_{1}(0^{-})=0^{-}$, the minus ($-$) superscript means that at the instant
$0^{-}$ clock $1$ is just about to receive the internal energy kick from its
escape mechanism.

\item \label{work2}We consider that the initial phases of the three clocks
are: $\psi_{3}(0^{-})=\psi_{3}^{0}>\psi_{2}(0^{-})=\psi_{2}^{0}>0^{-}=\psi
_{1}(0^{-})=\psi_{1}^{0}$.

\item \label{work4}The perturbation satisfies the relation $P\left(
x+Px\right)  \simeq Px$ in first order.
\end{enumerate}

To obtain the desired model, we need to proceed through 6 steps, starting from
the following initial conditions, that is the phase differences of all pairs
of clocks, and considering them at various points of the cycle of the
reference clock.

In the sequel $\psi_{i}^{j}$ denotes the phase of clock $O_{i}$ at the $j-th$
step.

Next the six steps follow. We show all the details of the calculation in the Appendix.

\begin{center}
{INITIAL CONDITIONS}
\end{center}

The phase difference between $O_{3}$ and $O_{1}$ is
\[
(CA)_{0}={\psi}_{3}^{0}-{\psi}_{1}^{0}={\psi}_{3}^{0},
\]
and the phase difference between $O_{1}$ and $O_{3}$ is symmetric, in the
sense that
\[
(AC)_{0}={\psi}_{1}^{0}-{\psi}_{3}^{0}=-{\psi}_{3}^{0}=-(CA)_{0}.
\]
The phase difference between $O_{2}$ and $O_{1}$ is
\[
(BA)_{0}={\psi}_{2}^{0}-{\psi}_{1}^{0}={\psi}_{2}^{0}%
\]
and the phase difference between $O_{1}$ and $O_{2}$ is
\[
(AB)_{0}={\psi}_{1}^{0}-{\psi}_{2}^{0}=-{\psi}_{2}^{0}=-(BA)_{0}.
\]
The phase difference between $O_{3}$ and $O_{2}$ is
\[
(CB)_{0}={\psi}_{3}^{0}-{\psi}_{2}^{0}%
\]
and the phase difference between $O_{2}$ and $O_{3}$ is
\[
(BC)_{0}={\psi}_{2}^{0}-{\psi}_{3}^{0}=-(CB)_{0}.
\]

\hspace{3cm}

\begin{center}
STEPS LEADING TO THE CONSTRUCTION OF THE MODEL
\end{center}

STEP 1: first impact. Interactions of $O_{1}$ on $O_{2}$ and of $O_{1}$ on
$O_{3}$, at $t=0$.

When the system in position A attains phase $0$ $(\operatorname{mod}2\pi)$ it
receives a sudden supply of energy, for short \textquotedblleft a
kick\textquotedblright, from its escape mechanism, this kick propagates in the
common support of the three clocks and reaches the other two clocks.

Now, the phase difference between $O_{3}$ and $O_{1}$ is corrected by the
perturbative value $P$:
\[
\left(  CA\right)  _{I}=\left(  CA\right)  _{0}+P\left(  \left(  CA\right)
_{0}\right)  ={\psi}_{3}^{0}+P\left(  {\psi}_{3}^{0}\right)  ==-(AC)_{I},
\]

where $(AC)_{I}$ is the phase difference between $O_{1}$ and $O_{3}$, since
$P$ must be an odd function of the mutual phase difference.

 The phase
difference between $O_{2}$ and $O_{1}$ is
\[
\left(  BA\right)  _{I}=\left(  BA\right)  _{0}+P\left(  \left(  BA\right)
_{0}\right)  ={\psi}_{2}^{0}+P\left(  {\psi}_{2}^{0}\right)  =-\left(
AB\right)  _{I}.
\]

The phase difference between $O_{3}$ and $O_{2}$ depends on $\left(
CA\right)  _{I}$ and $\left(  BA\right)  _{I}$ and it is
\[
\left(  CB\right)  _{I}={\psi}_{3}^{0}-{\psi}_{2}^{0}+P({\psi}_{3}%
^{0})-P({\psi}_{2}^{0})=-(CA)_{I}.
\]

STEP 2: first natural time shift. The next clock to arrive at $2\pi^{-}$, from
working hypothesis 3.2 (\ref{work2}), is the clock $O_{3}$ at vertex $C$. The
situation right before $O_{3}$ receives its kick of energy is when the phase
of this clock is $2\pi^{-}$.

At this point we have%
\[%
\begin{cases}
\psi_{3}^{2} & ={2\pi}^{-}\\
\psi_{1}^{2} & =2\pi-\left(  {\psi}_{3}^{0}+P({\psi}_{3}^{0})\right) \\
\psi_{2}^{2} & =2\pi+{\psi}_{2}^{0}-{\psi}_{3}^{0}+P({\psi}_{2}^{0})-P({\psi
}_{3}^{0}).
\end{cases}
\]

STEP 3: second impact. Clock $O_{3}$ receives its internal kick, at the
position $2\pi$.

Now, we have%
\[%
\begin{cases}
\psi_{3}^{3} & ={2\pi}\\
\psi_{1}^{3} & \simeq2\pi-{\psi}_{3}^{0}-2P\left(  {\psi}_{3}^{0}\right) \\
\psi_{2}^{3} & \simeq2\pi+{\psi}_{2}^{0}-{\psi}_{3}^{0}+P\left(  {\psi}%
_{2}^{0}\right)  -P\left(  {\psi}_{3}^{0}\right)  +P\left(  {\psi}_{2}%
^{0}-{\psi}_{3}^{0}\right) \\
&
\end{cases}
\]

STEP 4: second natural time shift. The next clock to arrive at $2\pi^{-}$,
from working hypothesis 3.2 (\ref{work2}), is the clock $O_{2}$ at vertex $B$.
The situation right before $O_{2}$ receives its kick of energy is when the
phase of this clock is $2\pi^{-}$.

Then we have
\[%
\begin{cases}
\psi_{2}^{4} & =2\pi^{-}\\
\psi_{1}^{4} & \simeq2\pi-{\psi}_{2}^{0}-P\left(  {\psi}_{2}^{0}\right)
-P\left(  {\psi}_{3}^{0}\right)  -P\left(  {\psi}_{2}^{0}-{\psi}_{3}%
^{0}\right) \\
\psi_{3}^{4} & \simeq2\pi-{\psi}_{2}^{0}+{\psi}_{3}^{0}-P\left(  {\psi}%
_{2}^{0}\right)  +P\left(  {\psi}_{3}^{0}\right)  -P\left(  {\psi}_{2}%
^{0}-{\psi}_{3}^{0}\right)  .
\end{cases}
\]

STEP 5: third impact. Clock $O_{2}$ receives its internal energy kick. It
reaches the position $2\pi$.

Then we have
\[%
\begin{cases}
\psi_{2}^{5} & ={2\pi}\\
\psi_{3}^{5} & \simeq2\pi-{\psi}_{2}^{0}+{\psi}_{3}^{0}-P\left(  {\psi}%
_{2}^{0}\right)  +P\left(  {\psi}_{3}^{0}\right)  -2P\left(  {\psi}_{2}%
^{0}-{\psi}_{3}^{0}\right) \\
\psi_{1}^{5} & \simeq2\pi-{\psi}_{2}^{0}-2P\left(  {\psi}_{2}^{0}\right)
-P\left(  {\psi}_{3}^{0}\right)  -P\left(  {\psi}_{2}^{0}-{\psi}_{3}%
^{0}\right)  .
\end{cases}
\]

STEP 6 (the final): third natural time shift. The next clock to arrive at
$2\pi^{-}$, from working hypothesis 3.2 (\ref{work2}), is the clock $O_{1}$ at
vertex $A$. The situation before $O_{1}$ receives its kick of energy is when
the phase of this clock is $2\pi^{-}$, i.e., the cycles is complete.

At this point we are able to describe what happens to the phases after a
complete cycle of the reference clock.

We have
\[%
\begin{cases}
\psi_{1}^{6} & ={2\pi}^{-}\\
\psi_{2}^{6} & \simeq2\pi+{\psi}_{2}^{0}+2P\left(  {\psi}_{2}^{0}\right)
+P\left(  {\psi}_{3}^{0}\right)  +P\left(  {\psi}_{2}^{0}-{\psi}_{3}%
^{0}\right) \\
\psi_{3}^{6} & \simeq2\pi+{\psi}_{3}^{0}+P({\psi}_{2}^{0})+2P({\psi}_{3}%
^{0})-P({\psi}_{2}^{0}-{\psi}_{3}^{0}).
\end{cases}
\]

Now, computing the phase differences after the first cycle of $O_{1}$, we obtain

\[%
\begin{cases}
(BA)_{I} & =-(BA)_{0}+2P((BA)_{0})+P((CA)_{0})+P((BA)_{0}-(CA)_{0})\\
(CA)_{I} & =((CA)_{0})+P((BA)_{0})+2P((CA)_{0})-P((BA)_{0}-(CA)_{0})\\
&
\end{cases}
\]

Hence, if we set $x=BA$ and $y=CA$, we obtain the system
\[
\left\{
\begin{array}
[c]{c}%
x_{1}=x_{0}+2P(x_{0})+P(y_{0})+P(x_{0}-y_{0})\\
y_{1}=x_{0}+P(x_{0})+2P({y}_{0})-P(x_{0}-y_{0}).
\end{array}
\right.
\]

\begin{center}
THE MODEL
\end{center}

By iterating the argument above, we get, for $n$ equal to the number of cycles
described by $O_{1}$, the discrete dynamical system:%
\[
\left\{
\begin{array}
[c]{c}%
x_{n+1}=x_{n}+2P(x_{n})+P(y_{n})+P(x_{n}-y_{n})\\
y_{n+1}=y_{n}+P(x_{n})+2P({y}_{n})-P(x_{n}-y_{n}).
\end{array}
\right.
\]
If we write
\[
\left\{
\begin{array}
[c]{c}%
\varepsilon\varphi\left(  x,y\right)  =2P(x)+P(y)+P(x-y)\\
\varepsilon\gamma\left(  x,y\right)  =P(x)+2P({y})+P(y-x),
\end{array}
\right.
\]
then we have
\[
\varphi\left(  x,y\right)  =\gamma\left(  y,x\right)  \text{,}%
\]
and the iteration is a perturbation of the identity as
\[
\left[
\begin{array}
[c]{c}%
x_{n+1}\\
y_{n+1}%
\end{array}
\right]  =\left[
\begin{array}
[c]{cc}%
1 & 0\\
0 & 1
\end{array}
\right]  \left[
\begin{array}
[c]{c}%
x_{n}\\
y_{n}%
\end{array}
\right]  +\varepsilon\left[
\begin{array}
[c]{c}%
\varphi(x_{n},y_{n})\\
\varphi(y_{n},x_{n})
\end{array}
\right]  ,
\]
that we can also write as
\begin{equation}
X_{n+1}=F(X_{n})=X_{n}+\varepsilon\Omega(X_{n}), \label{Model}%
\end{equation}
where
\[
X_{n+1}=\left[
\begin{array}
[c]{c}%
x_{n+1}\\
y_{n+1}%
\end{array}
\right]  ,
\]

\[%
\begin{array}
[c]{c}%
F(X_{n})
\end{array}
=\left[
\begin{array}
[c]{cc}%
1 & 0\\
0 & 1
\end{array}
\right]  \left[
\begin{array}
[c]{c}%
x_{n}\\
y_{n},
\end{array}
\right]
\]
and
\[%
\begin{array}
[c]{c}%
\Omega(X_{n})
\end{array}
=\left[
\begin{array}
[c]{c}%
\varphi(x_{n},y_{n})\\
\varphi(y_{n},x_{n})
\end{array}
\right]  .
\]

We now consider $P\left(  x\right)  =\varepsilon\sin x$, where $\varepsilon
=\frac{\alpha\mu}{8h^{2}}$ from hypothesis \ref{Perturb1}, explicitly,%
\begin{align*}
\varphi\left(  x,y\right)   &  =2\sin x+\sin y+\sin\left(  x-y\right) \\
\gamma\left(  x,y\right)   &  =\sin x+2\sin y+\sin\left(  y-x\right)  .\\
&
\end{align*}

\section{Analysis of the model}\label{sec4}

\subsection{Fixed points and local stability}

In this section, we analyse the model (\ref{Model}) obtained in the previous
section. In a nutshell, in this section, we see that the system is
differentiable and invertible in $S=\left[  0,2\pi\right]  \times\left[
0,2\pi\right]  $ when $\varepsilon>0$ is small. The perturbation map
$\Omega\left(  x,y\right)  $ is periodic in $R^{2}$. This implies that the
solution of the problem in the set $S$ is a dynamical system and not the usual
semi-dynamical system associated with discrete time. That will provide a
reasonable simple structure to the problem of the stability of fixed points
and will enable to derive global properties. Moreover, we prove that for small
$\varepsilon$ the set $S$ is invariant for the dynamics of $F$, meaning that
the two phase differences of oscillators $O_{2}$ and $O_{3}$ relative to
oscillator $O_{1}$ stay in the interval $\left[  0,2\pi\right[  $.

In particular, the map $\Omega$ has the zeros $\left(  \pi,\pi\right)  $,
$\left(  \frac{2}{3}\pi,\frac{4}{3}\pi\right)  $ and $\left(  \frac{4}{3}%
\pi,\frac{2}{3}\pi\right)  $ in the interior of the set $S=\left[
0,2\pi\right]  \times\left[  0,2\pi\right]  $, which are fixed points of the
model $F$. There are also four trivial fixed points, $\left(  0,0\right)  $,
$\left(  0,2\pi\right)  $, $\left(  2\pi,0\right)  $ and $\left(  2\pi
,2\pi\right)  $ at the corners of $S$, and the four fixed points $\left(
0,\pi\right)  $, $\left(  \pi,2\pi\right)  $, $\left(  2\pi,\pi\right)  $ and
$\left(  \pi,0\right)  $ on the edges of $S$.

We now compute the Jacobian matrix $J\left(  x,y\right)  $ to establish the
dynamical nature of the fixed points in the usual way.

We have

\begin{equation}
J\left(  x,y\right)  =\left[
\begin{array}
[c]{cc}%
1 & 0\\
0 & 1
\end{array}
\right]  +\varepsilon\left[
\begin{array}
[c]{cc}%
2\cos x+\cos\left(  x-y\right)  & -\cos\left(  x-y\right)  +\ \cos y\\
\cos x-\cos\left(  x-y\right)  & \cos\left(  x-y\right)  +2\cos y
\end{array}
\right]  .\label{Jacobmatrix}%
\end{equation}

We first consider the fixed points of $F$ in the interior of $S$. We start
with $\left(  \frac{2}{3}\pi,\frac{4}{3}\pi\right)  $ and $\left(  \frac{4}%
{3}\pi,\frac{2}{3}\pi\right)  $. The Jacobian is exactly the same%

\[
\left[
\begin{array}
[c]{cc}%
1-3\frac{\sqrt{3}}{2}\varepsilon & 0\\
0 & 1-3\frac{\sqrt{3}}{2}\varepsilon
\end{array}
\right]  ,
\]
meaning that those two points are locally asymptotically stable for
$\varepsilon$ sufficiently small.

The Jacobian matrix of $F$ at $\left(  \pi\text{,}\pi\right)  $ is
\[
\left[
\begin{array}
[c]{cc}%
1-\varepsilon & -2\varepsilon\\
-2\varepsilon & 1-\varepsilon
\end{array}
\right]  ,
\]
with eigenvalues $1-3\varepsilon$ and $1+\varepsilon$, which qualifies
$\left(  \pi\text{,}\pi\right)  $ as a saddle point. The stable manifold has
direction $\left(  1,1\right)  $, and the unstable manifold is tangent at
$\left(  \pi\text{,}\pi\right)  $ to the vector $\left(  -1,1\right)  $.

We now consider now the points placed at the vertexes of $S$. The Jacobian
matrix of $F$ at $\left(  0\text{,}0\right)  $, $\left(  0,2\pi\right)  $,
$\left(  2\pi,0\right)  $ and $\left(  2\pi,2\pi\right)  $ is, for all of
them, the following
\[
\left[
\begin{array}
[c]{cc}%
1+3\varepsilon & 0\\
0 & 1+3\varepsilon
\end{array}
\right]  ,
\]
which qualifies all the vertexes of $S$ as a repellers.

On the vertical edges of $S$ we have the fixed points $\left(  0,\pi\right)
$, and $\left(  2\pi,\pi\right)  $, at which the Jacobian matrix of $F$ is
\[
\left[
\begin{array}
[c]{cc}%
1+\varepsilon & 0\\
2\varepsilon & 1-3\varepsilon
\end{array}
\right]  ,
\]
which qualifies $\left(  0,\pi\right)  $, and $\left(  2\pi,\pi\right)  $ as
saddle points. The stable manifold has the direction of the $y$ axis and the
unstable manifold is tangent at $\left(  0\text{,}\pi\right)  $ and $\left(
2\pi,\pi\right)  $ to the vector $\left(  2,1\right)  $.

Finally, at the horizontal edges of $S$ we have the Jacobian matrix of $F$ at
$\left(  \pi,0\right)  $, and $\left(  \pi,2\pi\right)  $
\[
\left[
\begin{array}
[c]{cc}%
1-3\varepsilon & 2\varepsilon\\
0 & 1+\varepsilon
\end{array}
\right]  ,
\]
which qualifies $\left(  \pi,0\right)  $ and $\left(  \pi,2\pi\right)  $ again
as saddle points. The stable manifold is the direction of the $x$ axis and the
unstable manifold is tangent at $\left(  \pi,0\right)  $ and $\left(  \pi
,2\pi\right)  $ to the vector $\left(  1,2\right)  $.

The local analysis of the fixed points of $F$ reveals a very symmetric
picture. When $\varepsilon>0$ is small ($0<\varepsilon<\varepsilon_{0}%
=\frac{1}{9}$ is good enough), $F$ is a small perturbation of the identity,
$F(\partial S)=\partial S$, the restriction of $F$ to the boundary of $S$,
$\partial S$, is a bijection (see section 4 for more details), and the
Jacobian determinant of $F$ is never null in the interior of $S$. Therefore,
$F$ is invertible on S.

\subsection{Heteroclinic connections and invariant sets}

We focus our attention on the existence of invariant subsets of $S$ for the
dynamics of $F$. Additionally, we below prove that $S$ is itself an invariant
set for the dynamics of $F$.

Recall that an heteroclinic (sometimes called a heteroclinic connection, or
heteroclinic orbit) is a path in phase space which joins two different
equilibrium points. In the sequel, by sa-heteroclinic, rs-heteroclinic, and
ra-heteroclinic, we mean an heteroclinic orbit connecting a saddle point to an
attractor, an heteroclinic orbit connecting a repeller to a saddle point, and
an heteroclinic orbit connecting a repeller to an attractor, respectively.

Let $F$ be our model map in some set $T$ with two fixed points $p$ and $q$.
Let $M_{u}\left(  F,p\right)  $ and $M_{s}\left(  F,q\right)  $ be the stable
manifold and the unstable manifold (\cite{AlSaYo}: pages 78, 403) of the fixed
points $p$ and $q$, respectively. Then, if by $M$ we denote the heteroclinic
connecting $p$ and $q$, we have
\[
M\subseteq M_{s}\left(  F,p\right)  \cap M_{u}\left(  F,q\right)  .
\]
In particular, $M$ is invariant, the $\alpha$-limit and $\omega$-limit sets of
the points of $M$ is respectively $p$ and $q$ (\cite{AlSaYo}: page 331).

The other orbits, i.e., with initial conditions not in $M$, cannot cross the
heteroclinic connections when the map $F$ is invertible. In that case, it
would be violated the injectivity of the map. In the sequel, we study the
heteroclinics that connect saddle points to the attractors. Those
heteroclinics determine the nature of all the flow of the dynamical system in
the plane, due to the invertible nature of $F$.

\subsubsection{Vertical heteroclinics}

Consider the two vertical lateral edges of $S$, $s_{0}$ and $s_{1}$ that are
the sets $s_{k}=\left\{  \left(  x,y\right)  \in S:\left(  x=2k\pi\right)
\wedge0\leq y\leq2\pi\right\}  $, $k=0,1$. Consider the image of these
segments under $F$. If we write $F=(F_{1},F_{2})$, then
\[%
\begin{cases}
F_{1}\left(  2k\pi,y\right)  & =2k\pi+\varepsilon\sin y+\varepsilon\sin\left(
-y\right)  =2k\pi\\
F_{2}\left(  2k\pi,y\right)  & =y+2\varepsilon\sin y+\varepsilon\sin
y=y+3\varepsilon\sin y,
\end{cases}
\]
meaning that for $\varepsilon$ small enough the edges $s_{k}$, $k=0,1$, are
invariant, as already mentioned in section 2. Because of the initial
conditions, on each of the edges $s_{k}$, $k=0,1$, the dynamics is given by%
\[%
\begin{cases}
x_{n+1} & =2k\pi,\\
y_{n+1} & =y_{n}+3\varepsilon\sin y_{n}\text{,}%
\end{cases}
.
\]

\begin{figure}
[h]
\begin{center}
\includegraphics[
height=2.2528in,
width=2.2632in
]%
{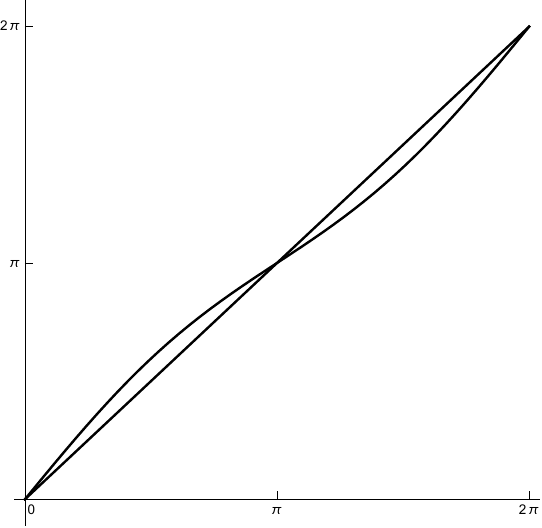}%
\caption{Graph of the map $g$, which is an homeomorphism in the interval
$\left[  0,2\pi\right]  $.}%
\label{homeo}%
\end{center}
\end{figure}

For $\varepsilon<\frac{1}{9}$, the map $g:\left[  0,2\pi\right]
\rightarrow\left[  0,2\pi\right]  $ defined by $g(t)=t+3\varepsilon\sin t$ is
a homeomorphism from the interval $\left[  0,2\pi\right]  $ into itself, as we
can see in Fig. \ref{homeo}. Moreover, since there is an attracting fixed
point of this map at $\pi$, the dynamics in the sets $s_{0}$ and $s_{1}$ can
be split in two subsets where the dynamics is again invariant, which is not
very important for our global discussion but establishes that the stable
manifolds of the saddle points $\left(  0,\pi\right)  $ and $\left(  2\pi
,\pi\right)  $ are, exactly and respectively, the sets $s_{0}$ and $s_{1}$

We have just shown that both $s_{0}$ and $s_{1}$ contain two heteroclinic
connections: in $s_{0}$, the line segment $s_{0}^{-}$ from $\left(
0,0\right)  $ to $\left(  0,\pi\right)  $ and $s_{0}^{+}$ from $\left(
0,2\pi\right)  $ to $\left(  0,\pi\right)  $; and in $s_{1}$, the line segment
$s_{0}^{-}$ from $\left(  2\pi,0\right)  $ to $\left(  2\pi,\pi\right)  $ and
$s_{0}^{+}$ from $\left(  2\pi,2\pi\right)  $ to $\left(  2\pi,\pi\right)  $.
The total number of vertical $rs$-heteroclinics is $4.$

\subsubsection{Horizontal heteroclinics}

Consider the two horizontal top and bottom edges of $S$, $r_{0}$ and $r_{1}$,
that are the sets $r_{k}=\left\{  \left(  x,y\right)  \in S:0\leq x\leq
2\pi\wedge\left(  y=2k\pi\right)  \right\}  $, $k=0,1$. Consider the image of
these segments under $F$. As before, if we write $F=(F_{1},F_{2})$, then
\[%
\begin{cases}
F_{1}\left(  x,2k\pi\right)  & =x+3\varepsilon\sin x\\
F_{2}\left(  x,2k\pi\right)  & =2k\pi,
\end{cases}
\]
meaning that, for $\varepsilon$ small enough, the edges $r_{k}$, $k=0,1$, are
invariant. Because of the initial conditions, on each of the edges $s_{k}$,
$k=0,1$, the dynamics is given by%
\[%
\begin{cases}
x_{n+1} & =x_{n}+3\varepsilon\sin x_{n}\\
y_{n+1} & =2k\pi\text{.}%
\end{cases}
\]
For $\varepsilon<\frac{1}{9}$, the map $g:\left[  0,2\pi\right]
\rightarrow\left[  0,2\pi\right]  $, defined by $g(t)=t+3\varepsilon\sin t$,
is the same occurred before, now involved in the dynamics in the invariant
edges $r_{0}$ and $r_{1}$. The stable manifolds of $\left(  \pi,0\right)  $
and $\left(  \pi,2\pi\right)  $ are again, respectively, the edges $r_{0}$ and
$r_{1}$.

Arguing as for $s_{0}$ and $s_{1}$, we have that both, $r_{0}$ and $r_{1}$,
contain two analogous heteroclinic connections.

We have just proved, in detail, that the boundary of $S$ is an invariant set.
More is true: each edge of $S$ is an invariant set.

Since the map $F$ is invertible, the initial conditions in the interior of
$S$, $S^{0}$, cannot cross the invariant boundary $\partial S=s_{0}\cup
s_{1}\cup r_{0}\cup r_{1}$, meaning that $S^{0}$ is an invariant set. This
means, in particular, that for equal clocks there will be no secular drift of
phase differences of the three clocks, the delays and advances are contained
in the set $S=\left[  0,2\pi\right]  \times\left[  0,2\pi\right]  $.

The total number of horizontal $rs$-heteroclinics is $4$. The total number of
$rs$-heteroclinics in the boundary of $S$ is $8$.

\subsubsection{Diagonal heteroclinics}

Finally, we now show that, $S^{o}$, the interior set of $S$, can be split
in two subsets, $S_{U}$ and $S_{D}$, $U$ for up and $D$ for down, where the
dynamics is again invariant. Consider now the set
\[
\Delta=\left\{  \left(  x,y\right)  \in S:y=x\text{, }x\in\left[
0,2\pi\right]  \right\}  ,
\]
the diagonal of $S$ connecting $\left(  0,0\right)  $ to $\left(  2\pi
,2\pi\right)  $. The image of a point of $\Delta$ by $F$ is now
\[%
\begin{cases}
F_{1}\left(  x,x\right)  & =x+3\varepsilon\sin x,\\
F_{2}\left(  x,x\right)  & =x+3\varepsilon\sin x.
\end{cases}
\]
Hence, the same homeomorphism $g$ as before appears again. We repeat the same
reasoning as before and deduce that $\Delta$ is invariant under $F$, and it
splits $S^{o}$ in two open sets: the triangle above it and the triangle below
it. Moreover, the stable manifold of the saddle point $\left(  \pi,\pi\right)
$ is the set $\Delta$.

This also proves the existence of two heteroclinics in $\Delta$, connecting
$\left(  0,0\right)  $ to $\left(  \pi,\pi\right)  $ and $\left(  2\pi
,2\pi\right)  $ to $\left(  \pi,\pi\right)  $, respectively. The total number
of $rs$-heteroclinics is now $10$, respectively $8$ on the edges and $2$ on the
main diagonal $\Delta$, all of them connecting repellers to saddles.

Consider now the other diagonal of $S$, i.e., the set
\[
\tilde{\Delta}=\left\{  \left(  x,y\right)  \in S:y=2\pi-x\text{, }x\in\left[
0,2\pi\right]  \right\}  .
\]
The image of a point of $\tilde{\Delta}$ under $F$ now is
\[%
\begin{cases}
F_{1}\left(  x,y\left(  x\right)  \right)  & =x+\varepsilon\sin x+\varepsilon
\sin2x,\\
F_{2}\left(  x,y\left(  x\right)  \right)  & =2\pi-\left(  x+\varepsilon\sin
x+\varepsilon\sin2x\right)  .
\end{cases}
\]
Hence, $\tilde{\Delta}$ is invariant.

\begin{figure}[h]
\begin{center}
\includegraphics[
height=2.3774in,
width=2.4457in
]%
{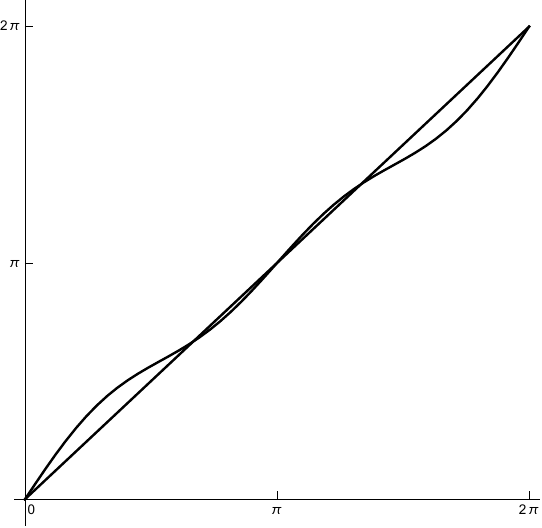}%
\caption{The homeomorphism $h_{1}$ with five fixed points.}%
\label{homeo2}%
\end{center}
\end{figure}

The map $h_{1}:\left[  0,2\pi\right]  \rightarrow\left[  0,2\pi\right]  $,
defined as $h_{1}(t)=t+\varepsilon\sin t+\varepsilon\sin2t$, is a
homeomorphism with $5$ fixed points from $\left[  0,2\pi\right]  $ to itself
(see Fig. \ref{homeo2}).

We repeat the same reasoning as before and deduce that the set $\tilde
{\Delta}$ splits the interior set $S^{o}$ again in two open sets: the triangle
above and the triangle below. So, now we have split $S^{0}$ in four small triangles.

There are four heteroclinic connections in $\tilde{\Delta}$, one connecting
the repeller $\left(  0,2\pi\right)  $ to the attractor $\left(  \frac{2\pi
}{3},\frac{4\pi}{3}\right)  $ (ra-heteroclinic), two $sa$-heteroclinics
connecting the saddle point $\left(  \pi,\pi\right)  $ to the attractors
$\left(  \frac{2\pi}{3},\frac{4\pi}{3}\right)  $ and $\left(  \frac{4\pi}%
{3},\frac{2\pi}{3}\right)  $, and, finally, the last heteroclinic on this
diagonal set is the one that connects the repeller $\left(  2\pi,0\right)  $
to the attractor $\left(  \frac{4\pi}{3},\frac{2\pi}{3}\right)  $
(ra-heteroclinic). The total number of sa-heteroclinics is now $2$.

\bigskip

We proceed with the same line of reasoning for the other $sa$-heteroclinics.

\begin{figure}
[h]
\begin{center}
\includegraphics[
height=2.3549in,
width=2.3929in
]%
{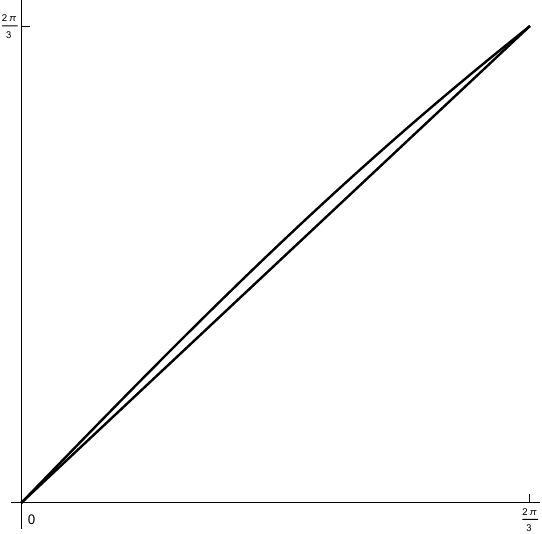}%
\caption{The homeomorphism $h_{1}$ with two fixed points, one repeller and the
other attractor.}%
\label{homeo3}%
\end{center}
\end{figure}

Consider now the set
\[
d_{1}=\left\{  \left(  x,y\right)  \in S:y=\pi+\frac{x}{2}\text{, }x\in\left[
0,\frac{2\pi}{3}\right]  \right\}
\]
and the map $F$ applied to the points of $d_{1}$:%
\[%
\begin{cases}
F_{1}\left(  x,y\left(  x\right)  \right)  & =x+2\varepsilon\sin
x-2\varepsilon\sin\left(  \frac{x}{2}\right)  ,\\
F_{2}\left(  x,y\left(  x\right)  \right)  & =\pi+\frac{1}{2}\left(
x+2\varepsilon\sin x-2\varepsilon\sin\left(  \frac{x}{2}\right)  \right)  .
\end{cases}
\]
The points of $d_{1}$ stay in $d_{1}$ under the action of $F$, proving that
this set also is invariant. The function $h_{2}:[0,2\frac{\pi}{3}%
]\rightarrow\lbrack0,2\frac{\pi}{3}]$, defined as $h_{2}(t)=t+2\varepsilon\sin
t-2\varepsilon\sin\left(  \frac{t}{2}\right)  $, is a homeomorphism, from
which we can readily see that the dynamics in $d_{1}$ is quite simple. The
graph of this homeomorphism can be seen in Fig. \ref{homeo3}. There is one
$sa$-heteroclinic from the saddle at $\left(  0,\pi\right)  $ to the attractor
$\left(  \frac{2\pi}{3},\frac{4\pi}{3}\right)  $. Actually, there is another
heteroclinic in the segment connecting the repeller $\left(  2\pi,2\pi\right)
$ to the attractor $\left(  \frac{2\pi}{3},\frac{4\pi}{3}\right)  $, but this
is not an sa-heteroclinic. Up to now, we have $3$ $sa$-heteroclinic connections.

Consider now the set $c_{1}=\left\{  \left(  x,y\right)  \in S:y=2x\text{,
}x\in\left[  \frac{2\pi}{3},\pi\right]  \right\}  $ and the map $F$ applied to
the points of $c_{1}$%
\begin{align*}
F_{1}\left(  x,y\right)   &  =x+\varepsilon\sin x+\varepsilon\sin2x,\\
F_{2}\left(  x,y\right)   &  =y+2\varepsilon\sin x+2\varepsilon\sin2x,
\end{align*}
the points of $c_{1}$ stay in $c_{1}$ under $F$, proving that this set is
invariant. Actually, the segment would be invariant if we extended $x$ to the
interval $\left[  0,\pi\right]  $, but we are not interested in heteroclinics
from repellers to attractors. Moreover, the dynamics is given by a restriction
of $h_{1}$ to the interval $\left[  \frac{2\pi}{3},\pi\right]  $. In this
interval there are only two fixed points, the attractor $\frac{2\pi}{3}$ and
the repeller $\pi$. This procedure adds one more $sa$-heteroclinic to the
global picture. So, we have found, up to now, $4$ $sa$-heteroclinics.

In $S_{D}$, we consider
\[
c_{2}=\left\{  \left(  x,y\right)  \in S:y=2\left(  x-\pi\right)  \text{,
}x\in\left[  \pi,\frac{4\pi}{3}\right]  \right\}
\]
and
\[
d_{2}=\left\{  \left(  x,y\right)  \in S:y=\frac{x}{2}\text{, }x\in\left[
\frac{4\pi}{3},2\pi\right]  \right\}  .
\]
Following exactly the same reasoning as before, we obtain two more
$sa$-heteroclinics, one connecting $\left(  \pi,0\right)  $ to the attractor
$\left(  \frac{4\pi}{3},\frac{2\pi}{3}\right)  $ and the other connecting
$\left(  2\pi,\pi\right)  $ to the same attractor.

\subsection{Phase portrait}

\begin{figure}[h]
\begin{center}
\includegraphics[
height=4.5in,
width=4.5in
]{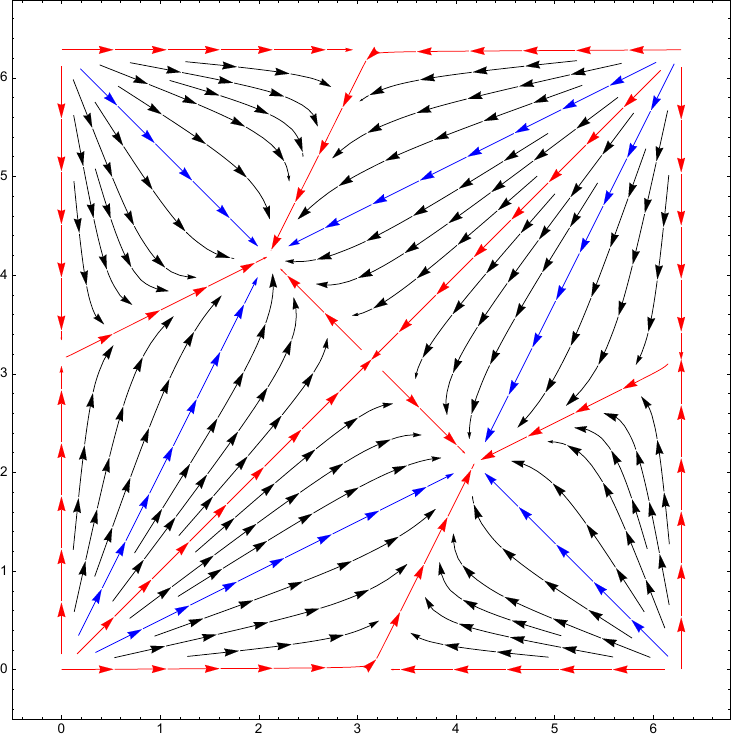}
\end{center}
\caption{Phase portrait of $F$, for small $\varepsilon$. In red the 16
saddle-node heteroclinic connections. For illustrative purposes, we represent
in blue some straight line invariant sets, actually heteroclinics, connecting
repellers to attractors. All the points in the interior of $S_{U}$ and $S_{D}$
belong to heteroclinics connections for $F$.}%
\label{Portrait}%
\end{figure}

The total number of $sa$-heteroclinics is $6$. All of them are straight
segments. The other $8$ $sa$-heteroclinics split the set $S$ in six invariant
sets as can be seen in Fig. \ref{Portrait}, where the red curves represent
saddle-node heteroclinics. The flow curves represented in the phase portrait.
Since the map $F$ is invertible, no orbit can cross either the red curves,
blue curves or black flow curves. There are only two attractors and the
dynamics, due to the invertible nature of the map $F$ and its large symmetry,
is relatively simple: in every invariant set in the plane, the restriction
maps are again homeomorphisms and the flow curves must follow, by continuity,
the heteroclinic connections on the outer boundaries of each invariant set.

Consequently, only the orbits on the outer edges and main diagonal, i.e., in
the set $s_{0}\cup s_{1}\cup r_{0}\cup r_{1}\cup d$ are not attracted to the
two attractors $\left(  \frac{2\pi}{3},\frac{4\pi}{3}\right)  $ and $\left(
\frac{4\pi}{3},\frac{2\pi}{3}\right)  $. The upper attractor $\left(
\frac{2\pi}{3},\frac{4\pi}{3}\right)  $ attracts the points in the open upper
triangle $S_{U}$ with converse results for the lower attractor $\left(
\frac{4\pi}{3},\frac{2\pi}{3}\right)  $ in $S_{D}$. The full picture can be
seen in Fig. \ref{Portrait}.

\subsection{Amplitude for three interacting oscillators in locked state}

\begin{figure}
[h]
\begin{center}
\includegraphics[
height=5.465625in,
width=4.6926in
]
{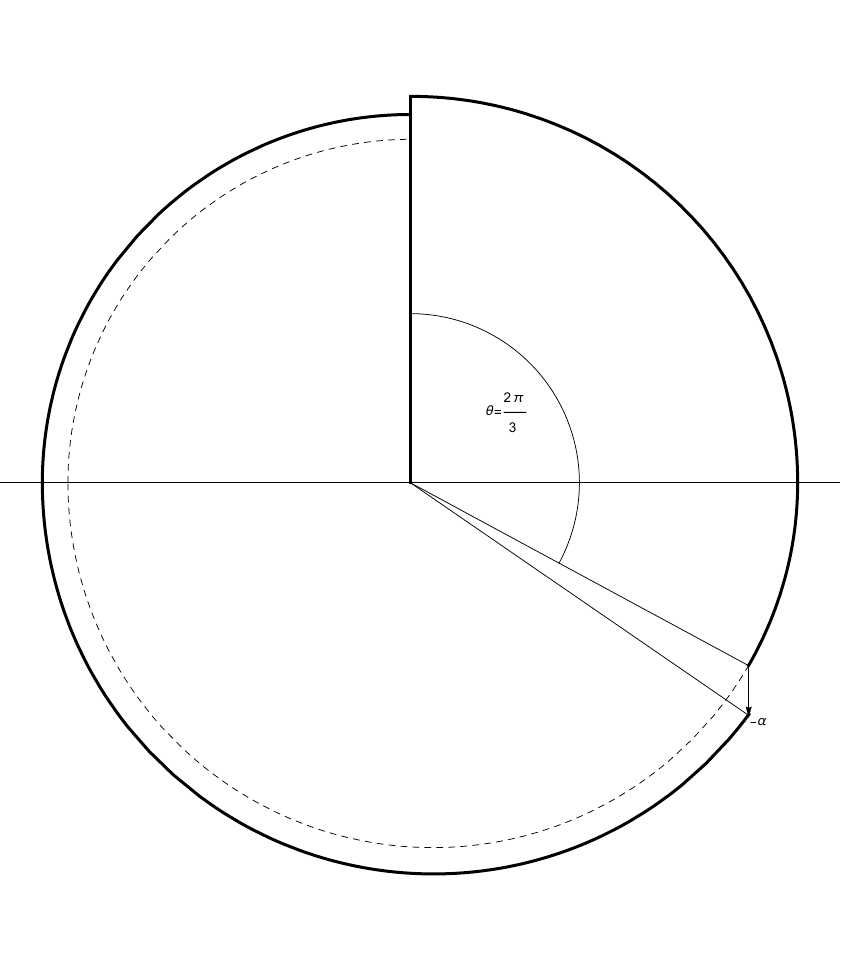}
\caption{When one of the oscillators receives a perturbation near the locked state at $\frac{2 \pi}{3}$ the radius of the limit cycle 
suffers an increase of ${\alpha}\sin(\frac{\pi}{6})=\frac{\alpha}{2}$, and the amplitude suffers an increase of the same value. Each 
oscillator suffers two impacts per cycle, therefore the total increase of the radius, i.e., the total increase of the amplitude, is $\alpha$. 
The effect of the perturbation is greatly exaggerated for the purpose of understandability.}
\label{Fig7}
\end{center}
\end{figure}

We focus again our
analysis on the velocity at $q=-\mu^{+}$ and $\dot{q}=v_{f}$. Since the mutual interaction affects each clock when they are at
relative phase differences of $\theta_{1}=\frac{2\pi}{3}$ and $\theta_{2}=\frac{4\pi}%
{3}$, as we can see in Fig. \ref{Fig7},  the equation for the Poincar\'{e} map of each generic oscillator
assumes now the form%
\begin{equation}
v_{n+1}=\sqrt{\left(  v_{n}-4\mu-\alpha\sin\left(  \theta_{1}-\frac{\pi}%
{2}\right)  -\alpha\sin\left(  \theta_{2}-\frac{\pi}{2}\right)  \right)
^{2}+h^{2}}\text{.}%
\end{equation}
Since%
\[
\sin\left(  \theta_{1}-\frac{\pi}{2}\right)  +\sin\left(  \theta_{2}-\frac
{\pi}{2}\right)  =1,
\]
this equation has again the same asymptotically stable\ fixed point
that we have seen in equation (\ref{perturb2}), which is again slightly
greater than the amplitude of the isolated clock%
\begin{equation}
v^{\ast}=\frac{\left(  4\mu+\alpha\right)  ^{2}+h^{2}}{2\left(  4\mu
+\alpha\right)  }\text{.}%
\end{equation}
The value of the maximum of the amplitude $v^{\ast}$ is the
same as in the case of two clocks.

The effect of the mutual perturbations between the clocks increases the amplitude of the oscillations, a phenomenon already 
observed in the experimental set up of \cite{OlMe}. 

\section{Conclusions and future work}\label{sec5}

In this paper we have proved that three oscillators, mutually interacting with
symmetric coupling, converge to a final symmetric locked state with mutual
phase differences of $\frac{2\pi}{3}$, this can happen in two different
settings, clockwise or counter-clockwise, depending on the initial conditions.

This very symmetrical final locked state induces us to consider the conjecture
that $n$ oscillators weakly interacting with all the others $n-1$ oscillators
will reach a final state with mutual phase differences of $\frac{2\pi}{n}$
clockwise or counter-clockwise distributed.

In future work, already in preparation, we shall discuss the same phenomenon
with slightly different natural angular frequencies $\omega_{1}$,
$\ \omega_{2}$ and $\omega_{3}$ and, in particular, the existence and form of
Arnold Tongues \cite{boyland1986,gilmore2011}.

As done for \cite{OlMe}, it would be interesting to check experimentally our
model, to see if the real world matches the theoretical predictions.

\subsection*{Disclosure of interest:} The authors report no conflict of interest.
\subsection*{Availability of data} Not applicable.

\bibliographystyle{plain}
\bibliography{BibloH2}

\section*{Appendix}

\subsection*{The steps in the construction of the model}

STEP 1: first impact. Interactions of $O_{1}$ on $O_{2}$ and of $O_{1}$ on
$O_{3}$, at $t=0$.

When the system in position A attains phase $0$ $(\operatorname{mod}2\pi)$ it
receives a sudden supply of energy, for short \textquotedblleft a
kick\textquotedblright, from its escape mechanism, this kick propagates in the
common support of the three clocks and reaches the other two clocks.

Now, the phase difference between $O_{3}$ and $O_{1}$ is corrected by the
perturbative value $P$:
\[
\left(  CA\right)  _{I}=\left(  CA\right)  _{0}+P\left(  \left(  CA\right)
_{0}\right)  ={\psi}_{3}^{0}+P\left(  {\psi}_{3}^{0}\right)  .
\]
The phase difference between $O_{1}$ and $O_{3}$ is
\[
\left(  AC\right)  _{I}=\left(  AC\right)  _{0}+P\left(  \left(  AC\right)
_{0}\right)  =-{\psi}_{3}^{0}+P\left(  -{\psi}_{3}^{0}\right)  =-(CA)_{I},
\]
since $P$ must be an odd function of the mutual phase difference.

The phase difference between $O_{2}$ and $O_{1}$ is
\[
\left(  BA\right)  _{I}=\left(  BA\right)  _{0}+P\left(  \left(  BA\right)
_{0}\right)  ={\psi}_{2}^{0}+P\left(  {\psi}_{2}^{0}\right)  ,
\]
and the symmetric phase difference between $O_{1}$ and $O_{2}$ is
\[
\left(  AB\right)  _{I}=\left(  AB\right)  _{0}+P\left(  \left(  AB\right)
_{0}\right)  =-{\psi}_{2}^{0}+P\left(  -{\psi}_{2}^{0}\right)  )=-\left(
BA\right)  _{I}.
\]

The phase difference between $O_{3}$ and $O_{2}$ depends on $\left(
CA\right)  _{I}$ and $\left(  BA\right)  _{I}$ and it is
\[
\left(  CB\right)  _{I}=\left(  CA\right)  _{I}-\left(  BA\right)  _{I}={\psi
}_{3}^{0}-{\psi}_{2}^{0}+P({\psi}_{3}^{0})-P({\psi}_{2}^{0})=-(CA)_{I},
\]

STEP 2: first natural time shift. The next clock to arrive at $2\pi^{-}$, from
working hypothesis 3.2 (\ref{work2}), is the clock $O_{3}$ at vertex $C$. The
situation right before $O_{3}$ receives its kick of energy is when the phase
of this clock is $2\pi^{-}$.

At this point we have%
\[%
\begin{cases}
\psi_{3}^{2} & ={2\pi}^{-}\\
\psi_{1}^{2} & =2\pi-(CA)_{I}=2\pi+(AC)_{I}=2\pi-\left(  {\psi}_{3}%
^{0}+P({\psi}_{3}^{0})\right) \\
\psi_{2}^{2} & =2\pi-\left(  CB\right)  _{I}=2\pi+\left(  BC\right)  _{I}%
=2\pi+{\psi}_{2}^{0}-{\psi}_{3}^{0}+P({\psi}_{2}^{0})-P({\psi}_{3}^{0}).
\end{cases}
\]

STEP 3: second impact. Clock $O_{3}$ receives its internal kick, at the
position $2\pi$.

Now, we have%
\[%
\begin{cases}
\psi_{3}^{3} & ={2\pi}\\
\psi_{1}^{3} & =\psi_{1}^{2}+P(\psi_{1}^{2})\\
& =2\pi-\left(  {\psi}_{3}^{0}+P({\psi}_{3}^{0})\right)  +P\left(
2\pi-\left(  {\psi}_{3}^{0}+P({\psi}_{3}^{0})\right)  \right) \\
& =2\pi-\left(  {\psi}_{3}^{0}+P({\psi}_{3}^{0})\right)  -P\left(  {\psi}%
_{3}^{0}+P({\psi}_{3}^{0})\right) \\
& \simeq2\pi-{\psi}_{3}^{0}-2P\left(  {\psi}_{3}^{0}\right) \\
\psi_{2}^{3} & =\psi_{2}^{2}+P(\psi_{2}^{2})\\
& =2\pi+{\psi}_{2}^{0}-{\psi}_{3}^{0}+P({\psi}_{2}^{0})-P({\psi}_{3}^{0})\\
& +P(2\pi+{\psi}_{2}^{0}-{\psi}_{3}^{0}+P({\psi}_{2}^{0})-P({\psi}_{3}^{0}))\\
& =2\pi+{\psi}_{2}^{0}-{\psi}_{3}^{0}+P({\psi}_{2}^{0})-P({\psi}_{3}^{0})\\
& +P({\psi}_{2}^{0}-{\psi}_{3}^{0}+P({\psi}_{2}^{0})-P({\psi}_{3}^{0}))\\
& \simeq2\pi+{\psi}_{2}^{0}-{\psi}_{3}^{0}+P\left(  {\psi}_{2}^{0}\right)
-P\left(  {\psi}_{3}^{0}\right)  +P\left(  {\psi}_{2}^{0}-{\psi}_{3}%
^{0}\right) \\
&
\end{cases}
\]

STEP 4: second natural time shift. The next clock to arrive at $2\pi^{-}$,
from working hypothesis 3.2 (\ref{work2}), is the clock $O_{2}$ at vertex $B$.
The situation right before $O_{2}$ receives its kick of energy is when the
phase of this clock is $2\pi^{-}$.

Then we have
\[%
\begin{cases}
\psi_{2}^{4} & =2\pi^{-}\\
\psi_{1}^{4} & =\psi_{1}^{3}+2\pi-\psi_{2}^{3}\\
& \simeq2\pi-{\psi}_{3}^{0}-2P\left(  {\psi}_{3}^{0}\right)  +2\pi\\
& -\left(  2\pi+{\psi}_{2}^{0}-{\psi}_{3}^{0}+P\left(  {\psi}_{2}^{0}\right)
-P\left(  {\psi}_{3}^{0}\right)  +P\left(  {\psi}_{2}^{0}-{\psi}_{3}%
^{0}\right)  \right) \\
& =2\pi-{\psi}_{2}^{0}-P\left(  {\psi}_{2}^{0}\right)  -P\left(  {\psi}%
_{3}^{0}\right)  -P\left(  {\psi}_{2}^{0}-{\psi}_{3}^{0}\right) \\
\psi_{3}^{4} & =\psi_{3}^{3}+2\pi-\psi_{2}^{3}\\
& \simeq2\pi+{2\pi}-\left(  2\pi+{\psi}_{2}^{0}-{\psi}_{3}^{0}+P\left(  {\psi
}_{2}^{0}\right)  -P\left(  {\psi}_{3}^{0}\right)  +P\left(  {\psi}_{2}%
^{0}-{\psi}_{3}^{0}\right)  \right) \\
& \simeq2\pi-{\psi}_{2}^{0}+{\psi}_{3}^{0}-P\left(  {\psi}_{2}^{0}\right)
+P\left(  {\psi}_{3}^{0}\right)  -P\left(  {\psi}_{2}^{0}-{\psi}_{3}%
^{0}\right)  .
\end{cases}
\]

STEP 5: third impact. Clock $O_{2}$ receives its internal energy kick. It
reaches the position $2\pi$.

Then we have
\[%
\begin{cases}
\psi_{2}^{5} & ={2\pi}\\
\psi_{3}^{5} & =\psi_{3}^{4}+P(\psi_{3}^{4})\\
& \simeq2\pi-{\psi}_{2}^{0}+{\psi}_{3}^{0}-P\left(  {\psi}_{2}^{0}\right)
+P\left(  {\psi}_{3}^{0}\right)  -P\left(  {\psi}_{2}^{0}-{\psi}_{3}%
^{0}\right) \\
& +P(2\pi-{\psi}_{2}^{0}+{\psi}_{3}^{0}-P\left(  {\psi}_{2}^{0}\right)
+P\left(  {\psi}_{3}^{0}\right)  -P\left(  {\psi}_{2}^{0}-{\psi}_{3}%
^{0}\right)  )\\
& \simeq2\pi-{\psi}_{2}^{0}+{\psi}_{3}^{0}-P\left(  {\psi}_{2}^{0}\right)
+P\left(  {\psi}_{3}^{0}\right)  -P\left(  {\psi}_{2}^{0}-{\psi}_{3}%
^{0}\right)  -P({\psi}_{2}^{0}-{\psi}_{3}^{0})\\
& =2\pi-{\psi}_{2}^{0}+{\psi}_{3}^{0}-P\left(  {\psi}_{2}^{0}\right)
+P\left(  {\psi}_{3}^{0}\right)  -2P\left(  {\psi}_{2}^{0}-{\psi}_{3}%
^{0}\right) \\
\psi_{1}^{5} & =\psi_{1}^{4}+P\left(  \psi_{1}^{4}\right) \\
& \simeq2\pi-{\psi}_{2}^{0}-P\left(  {\psi}_{2}^{0}\right)  -P\left(  {\psi
}_{3}^{0}\right)  -P\left(  {\psi}_{2}^{0}-{\psi}_{3}^{0}\right)  +\\
& P\left(  2\pi-{\psi}_{2}^{0}-P\left(  {\psi}_{2}^{0}\right)  -P\left(
{\psi}_{3}^{0}\right)  -P\left(  {\psi}_{2}^{0}-{\psi}_{3}^{0}\right)  \right)
\\
& \simeq2\pi-{\psi}_{2}^{0}-P\left(  {\psi}_{2}^{0}\right)  -P\left(  {\psi
}_{3}^{0}\right)  -P\left(  {\psi}_{2}^{0}-{\psi}_{3}^{0}\right)  -P({\psi
}_{2}^{0})\\
& =2\pi-{\psi}_{2}^{0}-2P\left(  {\psi}_{2}^{0}\right)  -P\left(  {\psi}%
_{3}^{0}\right)  -P\left(  {\psi}_{2}^{0}-{\psi}_{3}^{0}\right)  .
\end{cases}
\]

STEP 6 (the final): third natural time shift. The next clock to arrive at
$2\pi^{-}$, from working hypothesis 3.2 (\ref{work2}), is the clock $O_{1}$ at
vertex $A$. The situation before $O_{1}$ receives its kick of energy is when
the phase of this clock is $2\pi^{-}$, i.e., the cycles is complete.

At this point we are able to describe what happens to the phases after a
complete cycle of the reference clock.

We have
\[%
\begin{cases}
\psi_{1}^{6} & ={2\pi}^{-}\\
\psi_{2}^{6} & =\psi_{2}^{5}+2\pi-\psi_{1}^{5}\\
& \simeq2\pi+2\pi-\left(  2\pi-{\psi}_{2}^{0}-2P\left(  {\psi}_{2}^{0}\right)
-P\left(  {\psi}_{3}^{0}\right)  -P\left(  {\psi}_{2}^{0}-{\psi}_{3}%
^{0}\right)  \right) \\
& =2\pi+{\psi}_{2}^{0}+2P\left(  {\psi}_{2}^{0}\right)  +P\left(  {\psi}%
_{3}^{0}\right)  +P\left(  {\psi}_{2}^{0}-{\psi}_{3}^{0}\right)  ;\\
\psi_{3}^{6} & =\psi_{3}^{5}+2\pi-\psi_{1}^{5}\\
& \simeq2\pi-{\psi}_{2}^{0}+{\psi}_{3}^{0}-P\left(  {\psi}_{2}^{0}\right)
+P\left(  {\psi}_{3}^{0}\right)  -2P\left(  {\psi}_{2}^{0}-{\psi}_{3}%
^{0}\right)  +2\pi\\
& -\left(  2\pi-{\psi}_{2}^{0}-2P\left(  {\psi}_{2}^{0}\right)  -P\left(
{\psi}_{3}^{0}\right)  -P\left(  {\psi}_{2}^{0}-{\psi}_{3}^{0}\right)  \right)
\\
& =2\pi+{\psi}_{3}^{0}+P({\psi}_{2}^{0})+2P({\psi}_{3}^{0})-P({\psi}_{2}%
^{0}-{\psi}_{3}^{0}).
\end{cases}
\]

Now, we compute the phase differences after the first cycle of $O_{1}$.

We have
\begin{align*}
(BA)_{I}  &  =-(AB)_{I}=\psi_{2}^{6}-\psi_{1}^{6}\\
&  \simeq2\pi+{\psi}_{2}^{0}+2P\left(  {\psi}_{2}^{0}\right)  +P\left(  {\psi
}_{3}^{0}\right)  +P\left(  {\psi}_{2}^{0}-{\psi}_{3}^{0}\right)  -2\pi\\
&  ={\psi}_{2}^{0}+2P\left(  {\psi}_{2}^{0}\right)  +P\left(  {\psi}_{3}%
^{0}\right)  +P\left(  {\psi}_{2}^{0}-{\psi}_{3}^{0}\right) \\
&  =(BA)_{0}+2P((BA)_{0})+P((CA)_{0})+P((BA)_{0}-(CA)_{0})\\
&
\end{align*}
and%
\begin{align*}
&  (CA)_{I}\\
&  =-(AC)_{I}=\psi_{3}^{6}-\psi_{1}^{6}\\
&  =2\pi+{\psi}_{3}^{0}+P({\psi}_{2}^{0})+2P({\psi}_{3}^{0})-P({\psi}_{2}%
^{0}-{\psi}_{3}^{0})-2\pi\\
&  ={\psi}_{3}^{0}+P({\psi}_{2}^{0})+2P({\psi}_{3}^{0})-P({\psi}_{2}^{0}%
-{\psi}_{3}^{0})\\
&  =({(CA)}_{0})+P({(BA)}_{0})+2P({(CA)}_{0})-P((BA)_{0}-(CA)_{0})\\
&
\end{align*}
Hence, if we set $x=BA$ and $y=CA$, we obtain the system
\[
\left\{
\begin{array}
[c]{c}%
x_{1}=x_{0}+2P(x_{0})+P(y_{0})+P(x_{0}-y_{0})\\
y_{1}=x_{0}+P(x_{0})+2P({y}_{0})-P(x_{0}-y_{0}).
\end{array}
\right.
\]

\end{document}